\documentclass{article}

\usepackage[utf8]{inputenc}
\usepackage[small]{my-dgruyter}
\usepackage{microtype}

\usepackage{amsmath,amsfonts,amsthm,amssymb,eucal}

\usepackage{graphics} 
\usepackage{epsfig} 
\usepackage{graphicx}  
\usepackage{epstopdf}

\def\a{\alpha}

\def\LD{{^CD_{a+}^{\a}}}
\def\RD{{^CD_{b-}^{\a}}}

\def\RRD{{D_{b-}^{\a}}}

\newtheorem{theorem}{Theorem}

\newtheorem{example}{Example}
\theoremstyle{definition}

\theoremstyle{remark}

\begin{document}

\author[1]{Ricardo Almeida}

\author[1]{Delfim F. M. Torres}

\runningauthor{R. Almeida and D. F. M. Torres}

\affil[1]{Center for Research and Development in Mathematics and Applications (CIDMA),
Department of Mathematics, University of Aveiro,
3810-193 Aveiro, Portugal,
e-mail: \{ricardo.almeida, delfim\}@ua.pt}

\title{A survey on fractional variational calculus}

\runningtitle{A survey on fractional variational calculus}

\abstract{Main results and techniques of the fractional 
calculus of variations are surveyed. 
We consider variational problems 
containing Caputo derivatives and study them using both indirect 
and direct methods. In particular, we provide necessary optimality conditions 
of Euler--Lagrange type for the fundamental, higher-order, and isoperimetric problems,
and compute approximated solutions based on truncated Gr\"{u}nwald--Letnikov 
approximations of Caputo derivatives.}



\maketitle


\section{Introduction}

Calculus of Variations is a prolific field of mathematics,
whose founders were Leonhard Euler (1707--1783) and Joseph Louis Lagrange 
(1736--1813). The name of the field was given by Euler in his 
\emph{Elementa Calculi Variationum} and is related with the main
technique used by Lagrange to derive first order necessary optimality 
conditions. Its development is in intimate connection with various problems 
occurring in biology, chemistry, control theory, dynamics, 
economics, engineering, physics, etc. The current development 
of the calculus of variations is based 
on modern analysis techniques \cite{MR2191744,MR2191745}.

The Fractional Variational Calculus (FVC) is devoted to generalizations 
of the classic calculus of variations and its modern face,
the theory of optimal control, to the case in which derivatives and integrals 
are understood as fractional operators of arbitrary order. 
The subject has been developed in the last two decades
and provides nontrivial generalizations of the calculus of variations and optimal control, 
which open doors to new and interesting modern scientific problems.
Indeed, the combination of the calculus of variations and FC is not artificial, 
coming naturally from more accurate descriptions of physical phenomena, in particular 
to better describe non-conservative systems in mechanics, where the inclusion of non-conservatism 
in the theory is extremely important \cite{MR1401316,MR1438729}.

This chapter is organized as follows.
In Section~\ref{sec:main}, we formulate
the fundamental problem of the FVC.
Then our survey proceed in two parts:
indirect methods, in Section~\ref{sec:03},
and a direct method, in Section~\ref{sec:04}.
In an indirect method, the Euler--Lagrange equations, 
and possible other optimality conditions as well,
are used. Typically, the indirect approach 
leads to a boundary-value problem that should be solved to determine 
candidate optimal trajectories, called extremals. 
Direct methods, in contrast, begin by approximating the variational problem,
reformulating it as a standard nonlinear optimization problem \cite{MyID:407}.

Here we restrict ourselves to scalar problems of Caputo type.
However, the techniques presented are easily adapted to the vectorial
case and other types of FDs and FIs. For a survey on the FVC with generalized 
fractional operators see \cite{MR3221831}. For more on this fruitful area, 
we refer the readers to the four books on the subject
and references therein: for a solid introduction to the FVC, see \cite{MR2984893};
for computational aspects, begin with \cite{MR3443073}; 
for advanced mathematical methods, including the question of existence of solutions, 
see \cite{MR3331286}; for the variable order FVC and the state of the art see \cite{book4}.


\section{The fundamental problem}
\label{sec:main}

The calculus of variations deals with functional optimization problems 
involving an unknown function $x$ and its derivative $x'$. 
The general way of formulating the problem is as follows:
$$
\min \mathcal{J}(x)=\int_a^b L(t,x(t),x'(t))\,dt
$$
subject to the boundary conditions
\begin{equation}
\label{classicalBound}
x(a)=x_a \quad \mbox{and} \quad x(b)=x_b, 
\quad \mbox{with}\, x_a,x_b\in\mathbb{R}.
\end{equation}
One way to solve this problem consists to determine the solutions 
of the second order differential equation
\begin{equation}
\label{classicalEL}
\partial_2L (t,x(t),x'(t))=\frac{d}{dt}\partial_3L (t,x(t),x'(t)), \quad t\in[a,b],
\end{equation}
together with the boundary conditions \eqref{classicalBound}. 
Equation \eqref{classicalEL} was studied by Euler, and later by Lagrange, 
and is known as the Euler--Lagrange equation.

The fractional calculus of variations is a generalization of the 
ordinary variational calculus, where the integer-order derivative 
is replaced by a fractional derivative $D^\a x$:
$$
\min \mathcal{J}(x)=\int_a^b L(t,x(t),D^\a x(t))\,dt.
$$
As Riewe noted in \cite{MR1401316}, 
``traditional Lagrangian and Hamiltonian mechanics cannot 
be used with nonconservative forces such as friction''.
FVC allows it in an elegant way.

Since there are several definitions of fractional derivatives, 
one finds several works on the calculus of variations for these 
different operators. Here we consider the Caputo fractional derivative \cite{Caputo:1967}, 
but analogous results can be formulated for other types of fractional derivatives. 
Our choice of derivative is based on the following important fact: 
the Laplace transform of the Caputo fractional derivative depends on integer-order 
derivatives evaluated at the initial time $t=a$. Thus, we have an obvious physical 
interpretation for them. For other types of differentiation, this may not be so clear. 
For example, when dealing with the Riemann--Liouville fractional derivative, 
the initial conditions are fractional.

The main problem of the fractional calculus of variations, in the context 
of the Caputo fractional derivative, is stated as follows. Given a functional
\begin{equation}
\label{funct1}
\mathcal{J}(x)=\int_a^b L\left(t,x(t),\LD x(t)\right)\,dt, \quad x\in\Omega,
\end{equation}
where $\Omega$ is  a given class of functions, find a curve $x$ for 
which $\mathcal{J}$ attains a minimum value. We suppose that $L:[a,b]\times \mathbb{R}^2\to\mathbb{R}$ 
is continuously differentiable with respect to the second and third variables
and, for every $x\in\Omega$, the function $\LD x$ exists and is continuous 
on the interval $[a,b]$. On the set of admissible functions of the problem, 
the boundary conditions
\begin{equation}
\label{bound}
x(a)=x_a \quad \mbox{and} \quad x(b)=x_b, 
\quad \mbox{with}\, x_a,x_b\in\mathbb{R},
\end{equation}
may be imposed. We observe that a curve $x^\star\in\Omega$ is a (local) 
minimizer of functional \eqref{funct1} if there exists some positive real 
$\delta$ such that, whenever $x\in\Omega$ and $\|x^\star-x\|<\delta$, 
one has $\mathcal{J}(x^\star)-\mathcal{J}(x) \leq 0$.


\section{Indirect methods}
\label{sec:03}

The literature about the fractional calculus
of variations is very vast \cite{MR3443073,book4,MR3331286,MR2984893}.
The usual procedure in the calculus of variations, to find a candidate 
to minimizer, consists to take variations of the optimal curve $x^\star$. 
Let $\epsilon$ be a real number close to zero and $\eta$ a function chosen 
so that $x^\star + \epsilon \eta$ belongs to $\Omega$. For example, 
if the functional is defined in a class of functions with the boundary 
conditions \eqref{bound}, then the curve $\eta$ must satisfy the conditions 
$\eta(a)=0$ and $\eta(b)=0$. The curve $x^\star + \epsilon \eta$ 
is called a variation of the curve $x^\star$. For simplicity of notation, 
we will use the operator $[\, \cdot \, ]^\alpha$ defined by
$$
[x]^\alpha(t):=\left(t,x(t),\LD x(t)\right).
$$
Suppose that $x^\star$ minimizes $\mathcal{J}$, and consider 
a variation $x^\star + \epsilon \eta$. Depending on the existence 
of boundary conditions, some constraints over $\eta(a)$ and $\eta(b)$ 
may be imposed. Consider a function $j$, defined on a neighborhood of zero, 
and given by
$$
j(\epsilon):=\mathcal{J}(x^\star+\epsilon \eta)
=\int_a^b L(t,x^\star(t)+\epsilon \eta(t),\LD x^\star(t)+\epsilon\LD \eta(t))\,dt.
$$
Since $\mathcal{J}$ attains a minimum value at $x^\star$, then we conclude that 
$\epsilon=0$ is a minimizer of $j$, and thus $j'(0)=0$. Since
$$
j'(\epsilon)=\int_a^b \partial_2 L[x^\star+\epsilon \eta]^\alpha(t)\eta(t)
+ \partial_3 L[x^\star+\epsilon \eta]^\alpha(t)\LD \eta(t)\,dt,
$$
we conclude that
\begin{equation}
\label{auxi1}
\int_a^b \partial_2 L[x^\star]^\alpha(t)\eta(t)+ \partial_3 L[x^\star]^\alpha(t)\LD \eta(t)\,dt=0.
\end{equation}
Using now fractional integration by parts \cite{MR2218073}, we obtain that
\begin{multline}
\label{auxi2}
\int_a^b \partial_3 L[x^\star]^\alpha(t)\LD \eta(t)\,dt\\
=\int_a^b \RRD \left( \partial_3 L[x^\star]^\alpha(t)\right) \eta(t)\,dt 
+\left[{I_{b-}^{1-\a}}\left(\partial_3 L[x^\star]^\alpha(t)\right)\eta(t)\right]_a^b.
\end{multline}
Combining equations \eqref{auxi1} and \eqref{auxi2}, 
we prove that
\begin{multline}
\label{auxi3}
\int_a^b \left[ \partial_2 L[x^\star]^\alpha(t)
+\RRD \left( \partial_3 L[x^\star]^\alpha(t)\right)\right]^\alpha \eta(t)\,dt\\
+\left[{I_{b-}^{1-\a}}\left(\partial_3 L[x^\star]^\alpha(t)\right)\eta(t)\right]_a^b=0.
\end{multline}
Assume that the set of admissible functions have to satisfy
the boundary constraints \eqref{bound}. 
In this case, the variation function $\eta$ must satisfy the conditions 
$\eta(a)=0$ and $\eta(b)=0$. Thus, equation \eqref{auxi3} becomes
$$
\int_a^b \left[ \partial_2 L[x^\star]^\alpha(t)
+\RRD \left( \partial_3 L[x^\star]^\alpha(t)\right)\right] \eta(t)\,dt=0.
$$
Assume now that the function
$t\mapsto \RRD \left( \partial_3 L[x^\star]^\alpha(t)\right)$
is continuous on $[a,b]$. Since $\eta$ is an arbitrary function 
on the open interval $]a,b[$, we conclude, by the fundamental 
lemma of the calculus of variations \cite{MR0160139}, that
\begin{equation}
\label{auxi4}
\partial_2 L[x^\star]^\alpha(t)+\RRD \left( \partial_3 L[x^\star]^\alpha(t)\right)=0
\end{equation}
for all $t\in[a,b]$. Assume now that $x(a)$ and $x(b)$ are free. 
By the arbitrariness of $\eta$, we can first consider variations such that 
$\eta(a)=0$ and $\eta(b)=0$, and in this case we prove equation \eqref{auxi4}. 
Therefore, replacing  \eqref{auxi4} in \eqref{auxi3}, we deduce that
$$
\left[{I_{b-}^{1-\a}}\left(\partial_3 L[x^\star]^\alpha(t)\right)\eta(t)\right]_a^b\,dt=0.
$$
If $\eta(a)=0$ and $\eta(b)\not=0$, then
$$
{I_{b-}^{1-\a}}\left(\partial_3 L[x^\star]^\alpha(t)\right)=0 \quad \mbox{at} \quad t=b,
$$
and if  $\eta(a)\not=0$ and $\eta(b)=0$, then
$$
{I_{b-}^{1-\a}}\left(\partial_3 L[x^\star]^\alpha(t)\right)=0 \quad \mbox{at} \quad t=a.
$$
In conclusion, we proved the following result, which provides a necessary optimality
condition, known as the fractional Euler--Lagrange equation, 
and two subsidiary conditions, known as transversality conditions.

\begin{theorem} 
Let $x^\star$ be an admissible function and suppose that 
$\mathcal{J}$ attains a minimum value at $x^\star$.
If there exists and is continuous the function 
$t\mapsto \RRD \left(\partial_3 L[x^\star]^\alpha(t)\right)$, then
\begin{equation}
\label{ELeq}
\partial_2 L[x^\star]^\alpha(t)+\RRD \left(\partial_3 L[x^\star]^\alpha(t)\right)=0
\end{equation}
for every $t\in[a,b]$. Moreover, if $x(a)$ is free, then
$$
{I_{b-}^{1-\a}}\left(\partial_3 L[x^\star]^\alpha(t)\right)=0 \quad \mbox{at} \quad t=a
$$
and, if $x(b)$ is free, then
$$
{I_{b-}^{1-\a}}\left(\partial_3 L[x^\star]^\alpha(t)\right)=0 \quad \mbox{at} \quad t=b.
$$
\end{theorem}

As observed in \cite{MR2453739}, using the well--know relation 
$$
\RRD f(t) = \RD f(t) + \RRD f(b)
= \RD f(t) + \frac{f(b)}{\Gamma(1-\a) (b-t)^\a},
$$
between Riemann--Liouville and Caputo derivatives, one can easily
rewrite the the Euler--Lagrange equation \eqref{ELeq} in the equivalent form
\begin{equation*}
\partial_2 L[x^\star]^\alpha(t)+\RD \left(\partial_3 L[x^\star]^\alpha(t)\right)
+ \frac{\left.\partial_3 L[x^\star]^\alpha(t)\right|_{t = b}}{\Gamma(1-\a) (b-t)^\a}
=0.
\end{equation*}
Following similar ideas, we can study the case for several dependent variables  
$(x_1(t),\ldots,x_m(t))$, with $m\in\mathbb{N}$. In this case, we obtain $m$ 
fractional differential equations, one for each function $x_i$, $i=1,\ldots,m$.

\begin{theorem} 
Let $\mathcal{J}$ be the functional
$$
\mathcal{J}(x_1,\ldots,x_m):=\int_a^b L[x]^\alpha_m(t)\,dt
$$
with
$$
[x]^\alpha_m(t):=\left(t,x_1(t),
\ldots,x_m(t),{^CD_{a+}^{\a_1}} x_1(t),\ldots,{^CD_{a+}^{\a_m}} x_m(t)\right),
$$
where $\a_1,\ldots,\a_m\in \, ]0,1[$ and $m\in\mathbb N$, defined on the set of functions
$$
\left\{x\in C^1([a,b],\mathbb R^m)\,: x(a)=x_a \wedge x(b)=x_b\right\},
$$
where $x_a,x_b$ are two fixed vectors in $\mathbb{R}^m$.
Suppose also that $L:[a,b]\times \mathbb{R}^{2m}\to\mathbb{R}$ 
is continuously differentiable with respect to its $i$th variable, $i=2,\ldots,2m+1$.
Let $(x_1^\star,\ldots,x_m^\star)$ be a minimizer of $\mathcal{J}$, and suppose that 
the functions $t\mapsto {D_{b-}^{\a_i}}\left(\partial_{i+1+m} L[x^\star]^\alpha_m(t)\right)$ 
exist and are continuous on $[a,b]$, $i=1,\ldots,m$. Then,
$$
\partial_{i+1} L[x^\star]^\alpha_m(t)+{D_{b-}^{\a_i}}\left(\partial_{i+1+m} 
L[x^\star]^\alpha_m(t)\right)=0
$$
for all $i=1,\ldots,m$ and for all $t\in[a,b]$.
\end{theorem}

In the next result, the lower bound of the cost functional is a real $A$, 
where $A>a$, that is, $A$ is greater than the lower bound 
of the fractional derivative. For simplicity of presentation, we restrict
ourselves to the case $m = 1$.

\begin{theorem} 
Define the functional $A$ as
$$
\mathcal{J}(x):=\int_A^b L[x]^\alpha(t)\,dt,
$$
where $A>a$. If $\mathcal{J}$ achieves a minimum value at $x^\star$, 
and if the maps $t\mapsto {D_{A-}^{\a}}\left(\partial_3 L[x^\star]^\alpha(t)\right)$, 
for $t\in[a,A]$, and $t\mapsto \RRD \left(\partial_3 L[x^\star]^\alpha(t)\right)$, 
for $[a,b]$, exist and are continuous, then
$$
\RRD \left(\partial_3 L[x^\star]^\alpha(t)\right)-{D_{A-}^{\a}}\left(\partial_3 
L[x^\star]^\alpha(t)\right)=0, \quad t\in[a,A],
$$
and
$$
\partial_2 L[x^\star]^\alpha(t)+\RRD \left(\partial_3 L[x^\star]^\alpha(t)\right)=0, 
\quad t\in[A,b].
$$
Moreover, the following transversality conditions are fulfilled:
$$
\left\{
\begin{array}{l}
{I_{b-}^{1-\a}}\left(\partial_3 L[x^\star]^\alpha(t)\right)
-{I_{A-}^{1-\a}}\left(\partial_3 L[x^\star]^\alpha(t)\right)=0 
\quad \mbox{at $t=a$, if $x(a)$ is free;}\\
{I_{A-}^{1-\a}}\left(\partial_3 L[x^\star]^\alpha(t)\right)=0
\quad \mbox{at $t=A$, if $x(A)$ is free;}\\
{I_{b-}^{1-\a}}\left(\partial_3 L[x^\star]^\alpha(t)\right)=0 
\quad \mbox{at $t=b$, if $x(b)$ is free.}
\end{array}
\right.
$$
\end{theorem}

For our next result, besides the boundary conditions, an integral constraint 
is imposed on the set of admissible functions. Such type of problems are known as 
isoperimetric problems, and the first example goes back to Dido, Queen of Cartage 
in Africa. She was interested if finding the shape of a curve with fixed perimeter, 
of maximum possible area. As we know, the solution is given by a circle. 
The calculus of variations presents a solution to such problem, 
with an integral constraint of type
$$
\int_a^b \sqrt{1+(x'(t))^2}\,dt=\mbox{constant.}
$$
In our case, we consider that the new constraint depends 
also on a fractional derivative, and it is of form
\begin{equation}
\label{isope}
G(x):=\int_a^b M\left(t,x(t),\LD x(t)\right)\,dt=K,\quad K\in\mathbb{R},
\end{equation}
where $M:[a,b]\times\mathbb{R}^2\to\mathbb{R}$ is a constant, 
such that there exist and are continuous the functions 
$\partial_2 M$ and $\partial_3 M$.

\begin{theorem} 
Suppose that $\mathcal{J}$ \eqref{funct1}, subject to the constraints \eqref{bound} 
and \eqref{isope}, attains a minimum value at $x^\star$. 
If $x^\star$ is not a solution for
\begin{equation}
\label{aux8}
\partial_2 M[x]^\alpha(t)+\RRD \left(\partial_3 M[x]^\alpha(t)\right)=0,
\quad \forall t \in[a,b],
\end{equation}
and if there exist and are continuous the functions 
$t\mapsto \RRD \left(\partial_{3} L[x^\star]^\alpha(t)\right)$ and $t\mapsto \RRD
\left(\partial_{3} M[x^\star]^\alpha(t)\right)$ on $[a,b]$, then there exists 
$\lambda\in\mathbb R$ such that $x^\star$ satisfies
$$
\partial_2 F[x]^\alpha(t)+\RRD \left(\partial_3 F[x]^\alpha(t)\right)=0, 
\quad \forall t \in[a,b],
$$
where we define the function $F$ as $F:=L+\lambda M$.
\end{theorem}

The case when $x^\star$ is a solution of \eqref{aux8} can be easily included:

\begin{theorem} 
Suppose that $\mathcal{J}$ \eqref{funct1}, subject to the constraints \eqref{bound} and \eqref{isope}, 
attains a minimum value at $x^\star$. If there exist and are continuous the functions 
$t\mapsto \RRD \left(\partial_{3} L[x^\star]^\alpha(t)\right)$ and $t\mapsto \RRD
\left(\partial_{3} M[x^\star]^\alpha(t)\right)$ on $[a,b]$, then there exist 
$\lambda_0,\lambda\in\mathbb{R}$, not both zero, such that $x^\star$ satisfies the equation
$$
\partial_2 F[x]^\alpha(t)+\RRD \left(\partial_3 F[x]^\alpha(t)\right)=0, 
\quad \forall t \in[a,b],
$$
where we define the function $F$ as $F:=\lambda_0 L+\lambda M$.
\end{theorem}

Next, we present another constrained type problem, but now the restriction is given by
\begin{equation}
\label{isoperi2}g(t,x(t))=0, 
\quad \forall t\in[a,b],
\end{equation}
where $x=(x_1,x_2)$ is a vector and 
$g:[a,b]\times\mathbb{R}^2\to\mathbb{R}$ is differentiable 
with respect to $x_1$ and $x_2$. Also, we have the boundary conditions
\begin{equation}
\label{bound2} 
x(a)=x_a \quad \mbox{and} \quad x(b)=x_b, 
\quad x_a,x_b\in\mathbb{R}^2.
\end{equation}

\begin{theorem} 
Consider the functional
$$
\mathcal{J}(x)=\int_a^b L[x]^\alpha_2(t)\,dt,
$$
where 
$$
[x]^\alpha_2(t):=\left(t,x_1(t),x_2(t),\LD x_1(t),\LD x_2(t)\right),
$$
defined on $C^1[a,b]\times C^1[a,b]$, subject to the constraints 
\eqref{isoperi2} and \eqref{bound2}. If $\mathcal{J}$ attains an extremum 
at $x^\star=(x_1^\star,x_2^\star)$, the maps
$t \mapsto\RRD  \left(\partial_{i+3} L[x^\star]^\alpha_2(t)\right)$, 
$i=1,2$, are continuous, and
$\partial_3 g(t,x(t))\not=0$ for all $t\in[a,b]$,
then there exists a continuous function 
$\lambda:[a,b]\to\mathbb{R}$ such that
$$
\partial_{i+1}L [x^\star]^\alpha_2(t)+\RD\left( \partial_{i+3}
L[x^\star]^\alpha_2(t)\right)+\lambda \partial_{i+1} g (t,x(t))=0
$$
for all $t\in[a,b]$ and $i=1,2$.
\end{theorem}

Infinite horizon problems are an important field of research, that deal 
with phenomena that spread along time \cite{MR3046396,MR3031162,MR2966852}. 
In this case, the cost functional 
is evaluated in an infinite interval $[a,\infty[$:
\begin{equation}
\label{funct2}
\mathcal{J}(x):=\int_a^\infty L[x]^\alpha\,dt,
\end{equation}
defined on a set of functions with fixed initial conditions: $x(a)=x_a$.
Since we are dealing with improper integrals, some attention to what 
is a minimizer is needed. We say that $x^\star$ is a local minimizer 
for $\mathcal{J}$ as in \eqref{funct2} if there exists some $ \epsilon>0$ such that, 
for all  $x$, if $\|x^\star-x\|<\epsilon$, then
$$
\lim_{T\to\infty}\inf_{b \geq T}\int_{a}^{b}[L[x^\star]^\alpha(t)
-L[x]^\alpha(t)] \, dt 
\leq 0.
$$
Also, let us define the functions
\begin{align*}
A(\epsilon, b) &:= \int_a^{b} \frac{L[x^\star+\epsilon v]^\alpha(t)
- L[x^\star]^\alpha(t)}{\epsilon}dt;\\
V(\epsilon, T) &:= \inf_{b \geq T}\int_a^{b} \left[L[x^\star+\epsilon v]^\alpha(t)
-L[x^\star]^\alpha(t)\right]dt;\\
W(\epsilon)&:= \lim_{T\to\infty} V(\epsilon,T),
\end{align*}
where $v \in C^1[a,\infty[$ is a function and $\epsilon$ a real.

\begin{theorem} 
Let $x^\star$ be a local minimal for $\mathcal{J}$ 
as in \eqref{funct2}. Suppose that:
\begin{enumerate}
\item $\displaystyle \lim_{\epsilon \to 0} \frac{V(\epsilon, T) }{\epsilon}$ exists for all $T$;
\item $\displaystyle \lim_{T\to\infty}\frac{V(\epsilon, T) }{\epsilon}$ exists uniformly for all $\epsilon$;
\item For every $T> a$ and $\epsilon\not=0$, there exists a sequence 
$\left(A(\epsilon, b_n)\right)_{n \in \mathbb{N}}$ such that
$$
\lim_{n \to \infty} A(\epsilon, b_n)= \inf_{b\geq T} A(\epsilon, b)
$$
uniformly for $\epsilon$.
\end{enumerate}
If there exists and is continuous the function $t\mapsto \RRD \left(
\partial_3 L[x^\star]^\alpha(t)\right)$ on $[a,b]$, for all $b>a$, then
$$
\partial_2 L[x^\star]^\alpha(t)+\RRD \left(\partial_3 L[x^\star]^\alpha(t)\right)=0,
$$
for all $b> a$. Also, we have
$$
\lim_{T\to\infty}\inf_{b \geq T}{I_{b-}^{1-\a}}\left(\partial_3 
L[x^\star]^\alpha(t)\right)=0 \quad \mbox{at} \quad t=b.
$$
\end{theorem}

So far, we considered variational problems with order $\alpha\in]0,1[$. 
Now we proceed by extending them to functionals depending on higher-order derivatives. 
To that purpose, consider the functional
\begin{equation}
\label{funHig}
\mathcal{J}(x):=
\int_a^b L\left(t,x(t),{^CD_{a+}^{\a_1}} x(t),\ldots, {^CD_{a+}^{\a_m}} x(t)\right)\,dt,
\end{equation}
defined on $C^m[a,b]$, such that $x^{(i)}(a)$ and $x^{(i)}(b)$ are fixed reals 
for $i\in\{0,1,\ldots,m-1\}$. Here, $m$ is a positive integer, $\a_i\in(i-1,i)$, 
for all $i\in\{1,\ldots,m\}$, and $L:[a,b]\times\mathbb R^{m+1}\to\mathbb{R}$ 
is differentiable with respect to the $i$th variable, for $i\in\{2,3,\ldots,m+1\}$.

\begin{theorem} 
Let $x^\star$ be a minimizer of $\mathcal{J}$ \eqref{funHig}, and suppose that for all 
$i\in\{1,\ldots,m\}$, there exist and are continuous the functions $t\mapsto 
{D_{b-}^{\a_i}}\left(\partial_{i+2} L[x^\star]^\alpha_m(t)\right)$ on $[a,b]$. Then,
$$
\partial_2 L[x^\star]^\alpha_m(t)+\sum_{i=1}^m{D_{b-}^{\a_i}}\left(
\partial_{i+2} L[x^\star]^\alpha_m(t)\right)=0
$$
for all $t\in[a,b]$, where 
$[x^\star]^\alpha_m(t):=\left(t,x(t),{^CD_{a+}^{\a_1}} x^\star(t),
\ldots, {^CD_{a+}^{\a_m}} x^\star(t)\right)$.
\end{theorem}

After solving the presented Euler--Lagrange equations, we still need to verify 
if we are in presence of a minimizer of the functional or not. We recall that 
the Euler--Lagrange equation is just a necessary optimality condition of first order, and thus 
its solutions may not be a solution for the problem. For the question of existence
of solutions, we refer to \cite{MR3071538,MR3200762}.
One possible way to check if we have a candidate for minimizer or maximizer 
is to apply the Legendre condition, which is a second order 
necessary optimality condition \cite{MR3784384,MR3223604}.
Under suitable convexity of the Lagrangian, sufficient 
conditions for global minimizer hold \cite{MR2736825}.

For the Legendre condition, we have to assume that the Lagrange function $L$ is such that 
its second order partial derivatives $\partial^2_{ij}L$, with $i,j\in\{2,3\}$, exist and are continuous.

\begin{theorem} 
Suppose that $x^\star$ is a local minimum for $\mathcal{J}$ as in \eqref{funct1}. Then,
$$
\partial^2_{33}L[x^\star]^\alpha(t)\geq0, \quad \forall  t \in[a,b]. 
$$
\end{theorem}

For our next result, we recall that a function $L:[a,b]\times\mathbb{R}^2\to\mathbb{R}$ 
is convex with respect to the second and third variables if
$$
L(t,x+v,y+w)-L(t,x,y)\geq \partial_2 L(t,x,y)v+\partial_3 L(t,x,y)w
$$
for all $t\in[a,b]$ and $x,y,v,w\in\mathbb{R}$.

\begin{theorem}  
If $L$ is convex and if $x^\star$ satisfies the Euler--Lagrange equation \eqref{ELeq}, 
then $x^\star$ minimizes $\mathcal{J}$ \eqref{funct1} 
when restricted to the boundary conditions \eqref{bound}.
\end{theorem}


\section{Direct methods}
\label{sec:04}

There are several different numerical approaches and methods 
to solve fractional differential equations
\cite{MR3513174,MR3157930,MR3294734}. 
Different discretizations of fractional derivatives are possible, 
but many of them do not preserve the fundamental properties 
of the systems, such as stability \cite{MR2656354,MyID:397}. 
For numerical calculations, a simple and powerful method preserving stability 
is obtained from Gr\"{u}nwald--Letnikov discretizations.
For a detailed numerical treatment of fractional differential equations,
based on Gr\"unwald--Letnikov fractional derivatives, we refer
the interested reader to \cite{Ort,MR2824679}. Here we just mention
that both explicit and implicit methods are possible:
Theorem~5.1 of \cite{MR2824679} shows that the explicit and implicit 
Gr\"unwald--Letnikov methods are asymptotically stable; while
Theorem~5.2 of \cite{MR2824679} gives conditions on the step size
for the explicit method to be absolute stable, asserting that
the implicit method is always absolute stable,
without any step size restriction.
It is also worth to underline the good convergence properties
and error behavior of Gr\"{u}nwald--Letnikov methods \cite{MR2824679}.

We start this section by recalling the (left) Gr\"{u}nwald--Letnikov fractional 
derivative of a function $x$. Let $\alpha>0$ be a real. 
The Gr\"{u}nwald--Letnikov fractional derivative of order $\alpha$
is defined by
$$
{_a^{GL}D_x^\alpha}x(t):=\lim_{h\rightarrow 0^+} 
\frac{1}{h^\alpha}\sum_{k=0}^\infty(-1)^k\binom{\alpha}{k}x(t-kh),
$$
where $\binom{\alpha}{k}$ stands for the generalization of binomial coefficients 
to real numbers. As usual, we will adopt the notation
$$
(w_k^\alpha):= (-1)^k\binom{\alpha}{k}.
$$
This fractional derivative is particularly useful to approximate 
the Caputo derivative. The method is now briefly explained.

Given an interval $[a,b]$ and a fixed integer $N$, let $t_j:=a+ jh$, 
$j = 0,1,\ldots,N$ and $h>0$, be a partition of the interval $[a,b]$. Then,
$$
\LD x(t_j)= \frac{1}{h^\alpha}\sum_{k=0}^j(w_k^\alpha) x(t_{j-k})
-\frac{x(a)}{\Gamma(1-\alpha)}(t_j-a)^{-\alpha}+O(h).
$$
Thus, truncated Gr\"{u}nwald--Letnikov fractional derivatives are first-order 
approximations of the Caputo fractional derivatives. The approximation 
used for the Caputo derivative is then
$$
\LD x(t_j)\approx  \frac{1}{h^\alpha}\sum_{k=0}^j(w_k^\alpha) x(t_{j-k})
-\frac{x(a)}{\Gamma(1-\alpha)}(t_j-a)^{-\alpha}:= \tilde{D} x(t_j).
$$
The next step is to discretize the functional \eqref{funct1}. 
For simplicity, let $h=(b-a)/N$ and let us consider 
the grid $t_j=a+jh$, $j=0,1,\ldots, N$. Then,
\begin{equation}
\label{disFuncl}
\begin{array}{ll}
\mathcal{J}(x)&=\displaystyle \sum_{k=1}^N\int_{t_{k-1}}^{t_k} L(t,x(t),\LD x(t))\,dt\\
&\approx \displaystyle \sum_{k=1}^N h L(t_k,x(t_k),\LD x(t_k))\\
&\approx \displaystyle \sum_{k=1}^N h L(t_k,x(t_k),\tilde{D} x(t_k))\,dt.
\end{array}
\end{equation}
The right hand side of \eqref{disFuncl} can be regarded as a function 
$\Psi$ of $N-1$ unknowns:
$$
\Psi(x_1,x_2,\ldots,x_{N-1})=\sum_{k=1}^N h L(t_k,x(t_k),\tilde{D} x(t_k)).
$$
To find an extremum for $\Psi$, one has to solve the following system 
of algebraic equations:
$$
\frac{\partial \Psi}{\partial x_i}=0,\qquad i=1,\ldots,N-1.
$$
Suppose that $h$ goes to zero. The solution obtained by this method converges 
to a function $x^\star$. Then, $x^\star$ is a solution of 
the Euler--Lagrange equation \eqref{ELeq} (for details,
see \cite[Theorem~8.1]{MR3443073}).

\begin{example}
\label{ex01}
Consider the functional
$$
\mathcal{J}(x)=\int_0^{10}\left({^CD_{0+}^{0.5}} x(t)-\frac{2}{\Gamma(3/2)}t^{1.5}\right)^2\,dt
$$
subject to the boundary conditions $x(0)=0$ and $x(10)=100$. Since ${^CD_{0+}^{0.5}}t^2
= \frac{2}{\Gamma(3/2)}t^{1.5}$, and $\mathcal{J}$ is nonnegative, we conclude that the function 
$x(t)=t^2$ is a minimizer of the functional. If we discretize the functional, for different values 
of $n$, we obtain several numerical approximations of $x$. In Table~\ref{tab1} we present the error 
of each $n$, where the error is given by the maximum of the absolute value 
of the difference between $x$ and the numerical approximation.
\begin{table}[h]
\centering
\caption{Errors for the numerical solution of problem of Example~\ref{ex01}.} \label{tab1}
\begin{tabular}{lcccc}\hline
$n$ & 10 & 50 & 100 & 200 \\  \hline
Error  & 0.425189 & 0.107622 & 0.056370 & 0.029215\\ \hline
\end{tabular}
\end{table}
\end{example}


\begin{example}
\label{ex02}
For our second example, we consider a problem where we do not 
know its exact solution. Let
$$
\mathcal{J}(x)=\int_0^{1}\left(x(t)\left({^CD_{0+}^{0.5}} x(t)\right)^2
-\sin(x(t))\right)^2\,dt
$$
subject to the boundary conditions $x(0)=0$ and $x(1)=1$. 
In Figure~\ref{fig1}, we present the results for $n=100$.
\begin{figure}[ht]
\centering
\includegraphics[width=6cm]{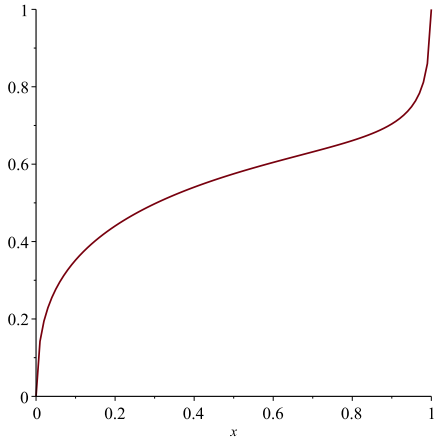}\\
\caption{Plot of the numerical solution of problem of Example~\ref{ex02}.}\label{fig1}
\end{figure}
\end{example}


\begin{acknowledgement}
the authors are grateful to the Research and Development 
Unit UID/MAT/04106/2013 (CIDMA) and to two anonymous referees
for several constructive comments.
\end{acknowledgement}



\end{document}